\documentstyle[12pt]{article}
\title{ Generalized Path Algebras and Pointed Hopf Algebras  }
\author{
Shouchuan Zhang $^{a,~b, ~c}$, \ \   Yao-Zhong Zhang $^b$, \ \
Xijing Guo $^a$ \
\\$a$. Department  of Mathematics, Hunan University\\ Changsha
410082, \
 P.R. China \\
$b$. Department of Mathematics, University of Queensland\\
Brisbane 4072, Australia\\
$c$. The Center of Combinatorial Mathematics, Nankai University\\
Tientsin, 300071, P.R.China\\ }
\date{}
\begin{document}
\newtheorem{Theorem}{\quad Theorem}[section]
\newtheorem{Proposition}[Theorem]{\quad Proposition}
\newtheorem{Definition}[Theorem]{\quad Definition}
\newtheorem{Corollary}[Theorem]{\quad Corollary}
\newtheorem{Lemma}[Theorem]{\quad Lemma}
\newtheorem{Example}[Theorem]{\quad Example}
\maketitle \addtocounter{section}{-1}

 \begin {abstract}

 Most of pointed Hopf
algebras of dimension $p^m$ with large coradical are shown to be
generalized path algebras. By the theory of generalized path
algebras it is  obtained that the representations, homological
dimensions and radicals of these Hopf algebras.  The relations
between the radicals of path algebras and connectivity of directed
graphs are given.

 \end {abstract}

\section {Introduction}

The concept of generalized path algebras was
introduced by F.U. Coelho and S.X. Liu in \cite {CL00}. Their  structures and
representations  were studied by the authors in \cite {ZZ03}.

Pointed Hopf algebras of dimension $p^m$ with large coradical were
classified (see \cite {DNR01}). In this paper, using the theory of
generalized path algebras we  give the representations,
homological dimensions and radicals of these Hopf algebras. In
Section 1, we  show that pointed Hopf algebra $H(C, n, c, c^*) $
is a generalized path algebra with relations and it is also a
smash product. The left global dimensions of Taft algebras are
infinity. The radicals and representations of $H(C, n, c, c^*)$
can be obtained by the theory of generalized path algebras in
\cite {ZZ03}. In Section 2, we give  the relations between
radicals of generalized matrix rings and $\Gamma$-rings. We obtain
the explicit formulas for generalized matrix ring $A$ and radical
properties $r = r_b,  r_l , r_j, r_n$:
$$r(A)=g.m.r(A) = \sum \{r (A_{ij}) \mid  i,  j \in I\}.$$
In Section 3, we  give the relations between the radicals of path
algebras and connectivity of directed graphs. That is, we obtain
that  every weak component of directed path D is a strong component
if and only if every unilateral component of $D$ is a strong
component if and only if the Jacobson radical of path algebra A(D)
is zero  if and only if  the Baer radical of path algebra A(D) is
zero.  In Section 4, we give some examples to establish that
generalized matrix algebras are applied in networks.

\vskip.1in
\noindent ${\bf Preliminaries }$
\vskip.1in
Let $k$ be a field.
We first recall the concepts of $\Gamma _I$-systems, generalized matrix rings
(algebras ) and generalized path algebras.
Let $I$ be a non-empty set. If for any $ i, j, l, s \in I, A_{ij}$
is an additive group and there exists a map
$\mu_{ijl}$ from $A_{ij}\times A_{jl}$ to $A_{il}$ (written
$\mu _{ijl} (x, y)=xy)$ such that the following conditions hold:

(i) $(x +y)z = xz +yz, \ \ \ w(x+y)= wx +wy;$

(ii) $w(xz)=(wx)z$,

\noindent for any $x, y \in A_{ij}$, $z\in A_{jl}, w\in A_{li}$, then the
set $\{A_{ij }\mid i,j \in I\}$ is a $\Gamma _I$ -system with index $I$.

Let $A$ be the external direct sum of $\{ A_{ij} \mid i, j\in I \}$. We define
the multiplication in $A$ as $$xy = \{ \sum _k x_{ik}y_{kj} \}$$
for any $x=\{x_{ij}\}, y=\{y_{ij}\}\in A$.
It is easy to check that $A$ is a ring (possibly without the unity element ).
We call $A$ a generalized matrix ring, or a gm ring in short,
written as $A=\sum \{A_{ij} \mid i, j\in I\}.$ For any non-empty subset $S $ of
$A$ and $i, j \in I$, set $S_{ij} = \{ a\in A_{ij}
\mid \hbox { there exists } x \in S \hbox { such that } x_{ij} = a \}$.
If $B$ is an ideal of $A$ and $B = \sum \{B_{ij} \mid i, j \in I\}$, then $B$ is called a gm ideal. If for any $i, j \in I$,
there exists $ 0 \not= e_{ii} \in A_{ii}$ such that $x_{ij}e_{jj} =
e_{ii}x _{ij } = x_{ij}$ for any $x_{ij} \in A_{ij}$,
then the set $\{e_{ii} \mid i \in I\}$ is called a generalized matrix unit of
$\Gamma _ I$-system $\{A_{ij} \mid i, j \in I\}$, or a generalized matrix  unit of  gm ring $ A = \sum \{A_{ij}
\mid i, j \in I\}$, or a gm unit in short. It is easy to show that
if $A$ has a gm unit $\{e_{ii} \mid i \in I\}$, then every ideal $B$ of $A$ is
a gm ideal. Indeed, for any
$x = \sum _{i, j \in I } x_{ij} \in B$ and $i_0, j_0 \in I$, since $e_{i_0i_0}
x e_{j_0 j_0} = x_{i_0j_0} \in B$,
we have $B _{i_0j_0} \subseteq B$. Furthermore, if $B$ is a gm ideal of $A$, then  $ \{A_{ij} / B_{ij} \mid i,
j\in I\}$ is a $\Gamma _I$-system and
$A/ B \cong \sum \{A_{ij} / B_{ij} \mid i, j\in I\}$ as  rings.

If for any $ i, j, l, s \in I, A_{ij}$
is a vector space over field $k$ and there exists a $k$-linear map
$\mu_{ijl}$ from $A_{ij}\otimes A_{jl}$ into $A_{il}$ (written
$\mu _{ijl} (x, y)=xy)$ such that $x(yz)=(xy)z$ for any $x\in
A_{ij}$, $y\in A_{jl}, z\in A_{ls}$, then the set $\{A_{ij } \mid i,j \in I \}$
is a $\Gamma _I$- system with index $I$ over field $k$. Similarly, we get an
algebra $A=\sum \{A_{ij} \mid i, j\in I\},$
called a generalized matrix algebra, or a gm algebra in short.

Assume that $D$ is a directed (or oriented) graph ($D$ is possibly an infinite
directed graph and also possibly not a simple graph) (or quiver ). Let $I=D_0$
denote the vertex set of $D$ and $D_1$ denote the set of arrows of
$D$. Let $\Omega $ be a generalized matrix algebra over field $k$ with
gm unit $\{e_{ii} \mid i \in I \}$, the Jacobson radical
$r(\Omega _{ii})$ of $\Omega _{ii}$ is zero and $\Omega _{ij} =0$ for
any $i \not= j \in I$.  The sequence
$x=a_{i_0} x_{i_0i_1} a_{i_1}x_{i_1i_2}a_{i_2} x_{i_2i_3} \cdots x_{i_{n-1}i_
{n}} a_{i_n}$ is called a generalized path (or $\Omega $-path) from $i_0$ to
$i_n$ via arrows $x_{i_0i_1}, x_{i_1i_2},x_{i_2i_3}, \cdots ,
x_{i_{n-1}i_{n}}$,
where $0\not= a_{i_p} \in \Omega _ {i_p i_p} $ for $p = 0, 1, 2, \cdots, n $. In this case, $n$ is called the length of $x$,
written $l(x).$
For two $\Omega $-paths
$x=a_{i_0} x_{i_0i_1}a_{i_1}x_{i_1i_2}a_{i_2}x_{i_2i_3} \cdots
x_{i_{n-1}i_{n}} a_{i_n}$ and
$y=b_{j_0}y_{j_0j_1} b_{j_1}y_{j_1j_2}b_{j_2}y_{j_2j_3} \cdots
y_{j_{m-1}j_{m}}b_{j_m}$ of $D$ with $i_n=j_0$, we define the
multiplication of $x$ and $y$ as
$$
xy = a_{i_0}x_{i_0i_1}a_{i_1}x_{i_1i_2}a_
{i_2} x_{i_2i_3} \cdots x_{i_{n-1}i_{n}}(a_{i_n}b_{j_0})y_{j_0j_1}
y_{j_1j_2} b_{j_1}y_{j_2j_3} \cdots y_{j_{m-1}j_{m}}b_{j_m}.  \ \ \ \ \ \ \ (*)
$$

\noindent For any $i, j \in I,$ let $A_{ij}'$ denote the vector
space over field $k$ with basis being all $\Omega $-paths  from
$i$ to $j$ with length $>0$. $B_{ij}$ is the sub-space spanned by
all elements of forms:
 $$ a_{i_0} x_{i_0i_1}a_{i_1}x_{i_1i_2}a_{i_2}\cdots x_{i_{s-1}i_s} (a_{i_s}^{(1)} + a_{i_s}^{(2)} + \cdots  +
 a_{i_s}^{(m)})x_{i_s i_{s+1}}
 \cdots
x_{i_{n-1}i_{n}} a_{i_n} $$
$$ - \sum _{l = 1}^ma_{i_0} x_{i_0i_1}a_{i_1}x_{i_1i_2}a_{i_2}x_{i_2i_3} \cdots x_{i_{s-1}i_s} a_{i_s}^{(l)} x_{i_s i_{s+1}} \cdots
x_{i_{n-1}i_{n}} a_{i_n},$$
where $i_0 =i, i_n = j,  a _{i_s}^{(l)} \in \Omega _{i_s i_s},  a _{i_p} \in \Omega _{i_p i_p}$,
 $x_{i_t i_{t+1}} $  is an arrow, $p = 0, 1, \cdots , n$, $t = 0, 1, \cdots , n-1$, $l = 0, 1, \cdots , m$,
$0 \leq s \leq n$,  $n$ and $m$ are
 natural
numbers. Let $A_{ij} = A_{ij}'/B_{ij}$ when $i\not= j$ and $A_{ii}
= (A_{ii}' +\Omega _{ii})/B_{ii}$, written $[\alpha ] = \alpha
+B_{ij}$ for any generalized path $\alpha $ from $i$ to $j$.
 We can get a $k$-linear map from $A_{ij}\otimes A_{jl}$ to $A_{il}$  induced by
$(*).$ We
write $a$ instead of  $[a]$ when $a\in \Omega$. In fact, $[\Omega _{ii}] \cong \Omega _{ii}$ as algebras for any $i\in I.$
 Notice that we write $e_{ii}x_{ij} = x_{ij}e_{jj} = x _{ij}$ for any arrow $x_{ij}$ from $i$ to $j$. It is clear that
$\{A_{ij} \mid i, j \in I\}$ is a $\Gamma _I$-system with gm unit
$\{e_{ii} \mid i \in I\}$. The gm algebra $\sum \{A_{ij}\mid i, j
\in I \}$ is called the generalized path algebra, or $\Omega
$-path algebra, written as $k (D, \Omega )$ (see, \cite [Chapter
3]{ARS95} and \cite {CL00}). Let $J$ denote the ideal generated by
all arrows in $D$ of $k(D, \Omega )$. If $\rho $ is a non-empty
subset of $k (D, \Omega )$ and the ideal $(\rho )$ generated  by
$\rho $ satisfies $J^t \subseteq (\rho ) \subseteq J^2 $, then
$k(D, \Omega )/(\rho )$ is called generalized path algebra with
relations. If $J^t \subseteq (\rho ) \subseteq J$, then $k(D,
\Omega )/(\rho )$ is called generalized path algebra with weak
relations. If $\Omega _{ii} = k e_{ii}$ for any $i\in I,$ then $k
(D, \Omega )$ is called a path algebra, written as $kD.$ For any
$i, j \in I $ and $A = k D$, if $u =\sum _{s= 1}^n  k_sp_s  \in
A_{ij}$, then the length  $l(u)$ is defined as the maximal length
$l (p_s)$ for $s= 1, 2, \cdots , n,$ where $p_1, p_2, \cdots p _s$
are different paths in $A_{ij}$. Furthermore, sometimes, we called
$k(D, \Omega )$  a generalized path algebra although the Jacobson
radical  $r(\Omega) \not=0.$

\section {Application in Hopf algebras}

Let $C$ be an abelian group and $G$ a group. Let $C^*$ denote the
character group of $C$,  ${\bf N}$ the set of natural numbers and
${\bf Z}^+ $ the set of positive integrals. Assume $c_i \in C$, $
c^*_i \in C^*, n = (n_1, n_2, \cdots , n_t)\in ({\bf Z}^+)^t$ and
$a = (a_1, a_2, \cdots, a_t) \in \{0, 1 \} ^t$; $ b_{ij} \in k$
for  $i, j =1, 2, \cdots , t.$ Throughout  this section, $D$
denotes the following quiver ( or directed graph ): vertex set
$D_0$ has only one element and arrow set $D_1 = \{X_1, X_2,
\cdots, X_t\}$. $\Omega = kC.$

Let $k_t = k \{ X_1, X_2, \cdot , X_t\}/(\rho )$ with
 $\rho = \{X_jX_i - c_j^*(c_i) X_iX_j \mid i, j = 1, 2, \cdots,
t \hbox { and } i\not= j\}$. That is, $k_t$ is the path algebra
 $k (D, \rho )$  with relation
$\rho $.

We first recall  $A_t(C, c, c^*, a, b)$ and $H(C, n, c, c^*, a,
b)$, which was defined in \cite [Definition 5.6.8 and  Definition
5.6.15] {DNR01}.

\begin {Definition} \label {3.3.1}
$A_t = A_t (C,  c, c^*, a, b)$ is the Hopf algebra generated by
the element $g\in C$ and $X_j , j = 1, 2, \cdots , t$ where

(i) the elements of $C$ are commuting grouplikes;

(ii)  $X_j$ is $(1, c_j)$-primitive;

(iii) $x_jg = c_j^* (g)gx_j;$

(iv) $X_jX_i = c^*_j(c_i) X_iX_j + b_{ij}(c_ic_j-1)$ for $ i, j
=1, 2, \cdots, t$ and $i \not= j;$

(v) $c^*_i(c_j)c^*_j(c_i) = 1$ for $j\not= i;$

(vi) If $b_{ij} \not= 0$ then $c^*_ic^*_j =1;$

(vii) If $c_ic_j =1,$ then $b_{ij} = 0.$

The antipode of $A_t$ is given  by $S(g) = g^{-1} $  for $g\in G$
and $S(X_j) = -c_j^{-1} X_j.$

Moreover,

(viii) $c_i^*(c_i)$  is a primitive $n_i$-th root of unity  for
any $i;$

(ix) If $a_i =1,$ then $(C^*_i)^{n_i} =1;$

(x) If $(c_i)^{n_i} =1 $ then $a_i =0;$

(xi) $b_{ij} = - c^*_i(c_j)b_{ji}$ for any $i, j.$

Let Hopf algebra  $H(C, n, c, c^*, a, b) = A_t/ J(a)$ where $J(a)$
is an ideal of $A_t$, generated by
$$ \{X_1^{n_1} - a_1 (c_1^{n_1}-1), X_2^{n_2} - a_2 (c_2^{n_2}-1),
\cdots, X_t^{n_t} - a_t (c_t^{n_t}-1) \}. $$

\end {Definition}

If $a=0$ and $b=0$, we denote  $A_t(C, c, c^*, a, b )$ and $H(C,
n, c, c^*, a, b )$ by $A_t(C, c, c^*)$ and $H(C, n, c, c^*)$,
respectively. If $C$ is a cyclic group generated by $c$ with order
$n$ and $t=1$ , then $H(C,n, c, c^* )$ is called a Taft algebra,
written as $H_{n^2} (\lambda )$, where $\lambda = c^*(c).$

Since we only consider the algebra structures of $A_t(C, c, c^*,
a, b )$ and $H(C, n, c, c^*, a, b)$ in this section, we use the
two signs ``$\cong $'' and ``='' denote isomorphism  and equation
as algebras, respectively.

\begin {Theorem}\label {3.3.2}
(i) $ A_t(C, c, c^*)=k_t \# kC .$

(ii) $A_t (C, c, c^*) = k(D, \Omega )/ ( \rho )$, with
 $\rho = \{ X_jX_i - c_j^*(c_i) X_iX_j, X_i h -c^*_i(h) h X_i \mid
 h \in C, i, j =1, 2, \cdots , t \hbox { and }$ {\ } $ i\not= j \}.$

(iii) $H(C, n, c, c^*) =  k(D, \rho) \# kC$ with $\rho = \{ X_jX_i
- c_j^*(c_i) X_iX_j, X_i^{n_i} \mid   i, j =1, 2, \cdots , t $ {\
\ \ } \hbox { and } $i\not= j \}.$

(iv) $H(C, n, c, c^*) = k (D, \Omega , \rho)$ with  $\rho = \{
X_i^{n_i}, X_jX_i - c_j^*(c_i) X_iX_j, X_i h -c^*_i(h) h X_i \mid
h\in C,   i, j =1, 2, \cdots , t \hbox { and } i\not= j \}.$

\end {Theorem}

{\bf Proof.} (i) For any $p = (p_1, p_2, \cdots , p_t) \in {\bf
Z}^t$, we denote  $X_1 ^{p_1}X_2 ^{p_2}\cdots X_t ^{p_t}$ by $X^p$
for convenience. Define $h\cdot X^p = C^*_1(h^{-1})^{p_2}
C^*_2(h^{-1})^{p_2} \cdots C^*_t(h^{-1})^{p_t}X^p  $ for any $h\in
C$. It is easy to check  this define a module structure of $k_t$
over $kC$ such that $k_t$ is a $kC$-module algebra.

Now we show that $k_t \# kC \cong A_t (C, c, c^* )$ as algebras.
Define  $\Phi : \ k_t \# kC \rightarrow  A_t (C, c, c^* ) $  via
sending $X^p \otimes g $ to $X^p g$ and $\Psi :   A_t (C, c, c^* )
  \rightarrow  k_t \# kG $ via sending $X^p g$ to $X^p \otimes g
$, for any $g, h \in C,$ $p = (p_1, p_2, \cdots , p _t), q = (q_1,
q_2, \cdots, q_t) \in {\bf N}^t$. See
\begin {eqnarray*}
( X^p \# g)(X^q\# h) &=& X^p (g \cdot X^q) \# gh\\
&=& (C^*_1(g^{-1})) ^{q_1}\cdots (C^*_t(g^{-1})) ^{q_t} X^p X^q \# gh, \hbox { \ \ \ and }\\
( X^p  g)(X^q h)  &=& (C^*_1(g^{-1})) ^{q_1}\cdots (C^*_t(g^{-1}))
^{q_t} X^p X^q  gh.
\end {eqnarray*}
Thus $\Phi $ is an algebra isomorphism. For convenience, we view
$A_t(C, c, c^*)$ as $k_t \# kC.$

(iii) Similarly, we can prove (iii).

(ii) and (iv) follow from the definition of  generalized path
algebras (see \cite {CL00}, \cite {ZZ03}). $\Box$




A representation of $(D, \Omega )$  is a set $ (V, f) =: \{ V, f_i
\mid V$ is an unitary $kC$-module, $f_i : V\rightarrow V$ is a
$k$-linear map, $ i= 1, 2, \cdots, n \}$. A morphism $h: (V, f)
\rightarrow (V', f')$ between two representations of $(D, \Omega
)$ is  a $k$-linear map $h : V \rightarrow V' $ such that  $h f_ i
= f_j 'h$ for $i , j = 1, 2, \cdots , t.$  Let {\cal R}ep $(D,
\Omega )$ denote the category of representations of $(D, \Omega
)$.

By \cite [Theorem 2.9] {ZZ03} we have

\begin {Corollary}\label {3.3.3} Let $\rho = \{ X_jX_i - c_j^*(c_i) X_iX_j, X_i h -
C^*_i (h) h X_i  \mid  i, j =1, 2, \cdots , t \hbox { and } i\not=
j \}.$ Then

(i) {\cal R}ep $k(D , \Omega , \rho )$ and ${}_{A_t (C, c, c^*)}
{\cal M}$ are equivalent.

(ii)  f.d.{\cal R}ep $(D, \Omega, \rho )$ and f.d.${}_{A_t(C, c,
c^*)} {\cal M} $ are equivalent. Here, f.d.{\cal R}ep $(D, \Omega
, \rho )$ and f.d.${}_{k(D , \Omega , \rho ) } {\cal M}$ denote
the full subcategories of finite dimensional objects in the
corresponding categories, respectively.

\end {Corollary}

Notice that although $ A_t (C, c, c^*) = k (D, \Omega ) /(\rho)$
is not a generalized path algebra with weak relations,  we may
use the conclusion in \cite [Theorem 2.9] {ZZ03}.
\begin {Corollary}\label {3.3.4} $\rho = \{ X_jX_i - c_j^*(c_i) X_iX_j,
x_i ^{n_i}, X_i h -c^*_i (h)hX_i  \mid i, j =1, 2, \cdots , t
\hbox { and } i\not= j \}.$
 Then

(i) {\cal R}ep $k(D , \Omega , \rho )$ and ${}_{H (C, n, c, c^*)}
{\cal M}$ are equivalent.

(ii)  f.d.{\cal R}ep $(D, \Omega, \rho )$ and f.d.${}_{H(C, n, c,
c^*)} {\cal M} $ are equivalent.
\end {Corollary}

Let $lgd (R)$ and $wd (R) $ denote the left global dimension and
weak dimension of algebra $R$, respectively.



\begin {Corollary}\label {3.3.5}   Let $k$ be a field  which
characteristic is  not divided by the order of finite group $C$.
Let $\rho = \{ X_jX_i - c_j^*(c_i) X_iX_j, X_i ^{n_i} \mid  i, j
=1, 2, \cdots , t \hbox { and } i\not= j \}.$ Then

(i) $r_b (A_t (C, c, c^*)) = 0.$

(ii) $lgd (A_t (C, c, c^*)) =   lgd  (k_t)$.

(iii) $ wd (A_t (C, c, c^*)) =   wd  ( k_t)$ .

(iv) $lgd (H(C,n,  c, c^*)) =   lgd  (k(D, \rho)) $.

(v) wd $(H (C, n,  c, c^*)) =  $ wd $  (k(D, \rho ))$.

\end {Corollary}

{\bf Proof.} (i) See that
\begin {eqnarray*}r_b (A_t (C, c, c^*)) &=& r_b (k_t \# kC) \hbox
{  \ by Theorem \ref {3.3.2}(i) }\\
&=& r_b (k_t) \# kC\hbox
{  \ by \cite [Theorem 2.6 ] {Zh98}  }\\
&=& 0 \hbox {\ \ \   since   } r_b (k_t) =0. \\
\end {eqnarray*}

(ii)--(v). By Theorem \ref {3.3.2} (i) (iii), $A_t (C, c, c^*)$
and $H(C, n, c, C^*)$ are smash products.  Using \cite [Theorem
2.2]{Zh01}, we can complete the proof. $\Box$

\begin {Corollary}\label {3.3.6} Let $k$ be a field  which
characteristic is  not divided by the order of finite group $C$.
Let $\rho = \{ X_jX_i - c_j^*(c_i) X_iX_j, X_i^{n_i}, X_i h-
C^*_i(h) h X_i \mid i, j =1, 2, \cdots , t \hbox { and } i\not= j
\}.$

(i) $r_j (H (C, n,  c, c^*)) = J/ (\rho),$ where $J$ is the ideal
of $H (C, n,  c, c^*)$, generated by $X_i's.$

(ii) $lgd (H_{n^2} (\lambda )) = \infty .$

\end {Corollary}
{\bf Proof.} (i) It follows from  Theorem  \ref {3.3.2} and \cite
[Lemma 3.6] {ZZ03}.

(ii) By Theorem \ref {3.3.2}, $H_{n^2} (\lambda )= (k[x]/ (x^n))\#
kC.$ By \cite [Corollary 2.2] {Zh01},
 $lgd (H_{n^2}(\lambda )) = lgd (k[x]/ (x^n)).$  It follows from
 The Hilbert syzygy Theorem $ lgd (k[x])=1$.
Let $R = k[x]$. If $ lgd (k[x]/(x^n)) < \infty , $ then
\begin {eqnarray*}
lgd ( R) &\geq &  \hbox {lgd } (R/ (x^n)) +1
\end {eqnarray*}  by \cite [Theorem 4.8] {Ch99} since $x^n$ is an
elelment in the center of $R$.  Consequently, $lgd (R/(x^n)) =0$.
This implies $R/(x^n)$ is a von Neumann regular ring, so $r_j
(R/(x^n))=0$. Obviously,  $(x)/(x^n)$ is a non-zero nilpotent
ideal. We get a contradiction. Therefore, lgd $(R/(x^n)) =
\infty.$ $\Box$

If we replace the abelian group $C$ by group $G$ in Definition
\ref {3.3.1}, we obtain two algebra $A_t(G, c, c^*, a, b)$ and
$H(G, n, c, c^*, a, b),$ where $c_i \in Z(G)$ for $i = 1, 2,
\cdots , t.$ It is clear that Theorem \ref {3.3.2} and Corollary
\ref {3.3.3} -- \ref {3.3.6} hold for $A_t(G, c, c^*, a, b)$ and
$H(G, n, c, c^*, a, b)$ since we only consider their algebra
structures.

\begin {Theorem}\label {3.3.7} Let $k$ be an algebraically closed  field of characteristic zero
and  $H$  a pointed Hopf algebra with dimension $p^m$, $p$ prime.
If $m \leq 3$
 or the dimension of the coradical of $H$ is more than $ p^{m-2}$, then $H$ is one of the following:

(i) a group algebra.;

(ii)   $H(C, n, c, c^*).$

(iii) $H(C, n, c, c^*, a, b).$

(iv)   $H(G, n, c, c^*)$ with $t=1.$

(v) $H(G, n, c, c^*, a, b)$ with $t=1.$

\end {Theorem}

{Proof.} We assume $H$ is not a group algebra.
If the dimension of the coradical of $H$ is more than $ p^{m-2}$,
then $G(H)=G$ is a group of order $p^{m-1}$. By \cite [Theorem
7.8.2] {DNR01}, $H= H(G, n,  c, c^*)$ with  $t= 1, n =p, c=g\in
Z(G)$, or $H= H(G, n,  c, c^*, a, 0)$ with $t= 1, n =p, c=g\in
Z(G)$ and $a=1$.

If $H$ is a Hopf algebra of dimension $p^3$ with $G(H)= (g)$ a
cyclic group of order $p$, then it follows from  \cite [Theorem
7.9.6] {DNR01} that $H\cong H(C, n, c, c^* )$ with $C=(g)$ of
order $p$, $t =2$, $n= (p, p), c= (g, g^i), c^* = (c^*,
(c^*)^{-i})$ for $1 \leq i \leq  p-1$, or
 $H\cong H(C, n, c, c^*, a, b )$ with $C= (g)$ of order $p$, $t =2$, $n= (p, p), c= (g, g),
c^* = (c^*, (c^*)^{-1}), a=0, b_{12}= 1. $ $\Box$

\begin {Example}\label {3.3.8}
Recall the duality theorem (see \cite [Corollary 6.5.6 and Theorem
6.5.11] {DNR01}) for co-Frobenius Hopf algebra $H$:
$$(R \# H^{* rat }) \# H\cong M_{H}^f (R) \hbox  { and } \\
(R \# H) \# H ^{* rat }\cong M_{H}^f (R) \ \ \ \ \  (as~ algebras
),$$ where $M_{H}^f (R) = \sum \{ B_{ij} \mid i, j \in I\}$ is a
gm algebra
 and $I$ is the basis of $H$ with $B_{ij} = R$ for $i, j \in I$. Note for
 the Baer radical, Levitzki radical, Jacobson radical and von
Neumann regular radical
  $r(M_H^f(R)) = M_H^f( r (R))$ of $M_H^f (R)$.  Consequently,
 $r((R \# H^{* rat }) \# H ) \cong  M_H^f( r (R))$  and      $r((R \# H) \# H ^{* rat }) \cong
  M_H^f( r (R)).$  In particular, the Heseberg algebra $H \# H^{* rat } {\ \ \ \
  }$ \ $
  \cong M_H^f( k) $  for
infinite co-Frobenius Hopf algebra H. Therefore
$$
  r (H \# H^{* rat }) \cong r(M_H^f(k)) = M_H^f( r( k)) = \left \{ \begin {array} {ll} 0
  & \ \ \hbox {when } r= r _b, r_l, r _j\\
  M_H^f(  k)    & \ \ \hbox {when } r= r _n.\\
  \end {array} \right. $$
\end {Example}

\section { The radicals of generalized matrix rings }

Since every generalized path algebra is a generalized matrix ring,
we study the radicals of generalized matrix rings in this section.

Let $x E(i,j)$ denote the generalized matrix
having a lone $x$ as its $(i, j)$-entry and all other entries are
zero. If $B$ is a non empty subset of generalized matrix ring $A$
and $s, t \in I$,  we call the set $\{x  \in A_{st} \mid  \hbox {
there exists }
 y \in  B  \hbox { such that  } y_{st} = x  \}$ the projection on $(s, t)$ of $B$, written as $B_{st}$.

Let $g.m.r(A)$ denote the maximal gm ideal of $r(A)$ for a radical
property $r$ of rings (see \cite {Zh93}). Let $r_b, r_l, r_k, r_j,
r_n $  denote the Baer radical, Levitzki radical, nil radical,
 Jacobson radical and von Neumann regular radical of rings and $\Gamma$ -rings, respectively. Let $r(A_{ij})$ denote
$r$  radical of $A_{ji}$-ring $A_{ij}$ for any $i, j \in I.$

Now we study the von Neumann radical $r_n (A)$ of generalized
matrix ring $A=\sum \{A_{ij} \mid i, j\in I \}.$
\begin {Definition}\label {1.1}

If for all $s, t \in I$, there exists  $0\not= d_{st}\in A_{st}$
such that $x_{is}d_{st}\not=0$ and $ d_{st} y_{tj} \not=0 $ for
any $i, j \in I, x_{is} \in A_{is}, y _{tj} \in A_{tj}$ , then we
say that $A$ has a left gm non-zero divisor.
\end {Definition}
Similarly, we can define the right gm non-zero divisor of $A$.

\begin {Lemma} \label {1.2}
(i)  If $B$ is an ideal of $A$, then $\bar B = \sum \{B_{ij}\mid
i,  j \in I\}$ is the gm ideal generated by $B$ in $A$.

(ii)  If $D$ is a gm ideal of $A$ and $D \subseteq \sum
\{r_n(A_{ij} )\mid  i,  j\in I\}$, then D is an $r_n$ -gm ideal of
$A$.

(iii)  Let $B_{st}$ be an $r_n$-ideal of $A_{ts}$-ring $A_{st}$
and $D_{ij}=A_{is}B_{st}A_{tj}$ for any $i, j\in I$. If $A$ has
 left and right gm non-zero divisors, then $D$ is an $ r_n$-gm
ideal of $A$.

(iv)  If  $A$ has
 left and right gm non-zero divisors and g.m. $r_n(A)=0$, then $r_n(A_{ij})=0$ for any $i, j\in I$.
\end {Lemma}

{\bf Proof.} (i) It is trivial.

(ii) For any $x\in D$, there exists a finite subset $J$ of $I$
such that $x_{ij} = 0 $ for any $i, j \notin J$. Without lost the
generality,  we can assume that $J = \{ 1, 2, \cdots , n  \}$ and
$J' = \{ 1, 2, \cdots , n, n+1 \}\subseteq I.$  Let $J'\times J' =
\{(u, v) \mid u, v = 1, 2, \cdots , n+1\}$ with the dictionary
order. We now show that there exist two sequences
 $\{ y _{t_2 t_1} \in A_{t_2 t_1} \mid (t_1,t_2)\in J'\times J' \}$
and   $\{x^{(t)} \in D  \mid (t_1,t_2) \in J'\times J' \}$ with $x
^{(1, 1)} = x$ and
\begin {eqnarray} \label {1.2e}
x^{(t +1)} = x^{(t)} - x^{(t)}(y _{(t_2t_1)} E(t_2,
t_1))x^{(t)}\end {eqnarray}
 such that
 $x ^{(t)}_s = 0$ for any  $ s, t \in J'\times J' $ with $s \prec t $  by induction. Since $x^{(11)}_{11}
 = x_{11}$ is a von Neumann regular  element, there exists $y_{11} \in A_{11}$ such that
$x _{11}= x_{11}y _{11} x_{11}$. See that $x^{(12) } _{11} =
x^{(11)}_{11}- x^{(11)} _{11} y_{11} x^{(11)}_{11}= 0$. For
$t=(t_1, t_2) \in B$, we assume that there exists $ y _{s_2 s_1}
\in A_{s_2 s_1}$ and   $x^{(t)}_{s} = 0$ for any $s = (s_1, s_2)
\prec (t_1, t_2)$.
  Since $x^{(t_1t_2)}_{t_1t_2}$
 is a von Neumann regular  element, there exists $y_{t_2t_1} \in A_{t_2t_1}$ such that
$x _{t_1t_2}^{(t_1t_2)}= x_{t_1t_2}^{(t_1 t_2)}y_{t_2t_1}
x_{t_1t_2}^{(t_1 t_2)}$. By (\ref {1.2e}), we have $x^{(t
+1)}_{t_1 t_2} = 0$.  For $s = (s_1,s_2)\prec  t = (t_1,t_2)$, we
have  either $s_1=t_1$, $s_2 < t_2$ or $s_1 <t_1$. This implies
that $(t_1, s_2 )\prec (t_1, t_2) $ or $(s_1, t_2 )\prec (t_1,
t_2) $. Thus $x_{s_1s_2} ^{(t +1)} = 0$ by (\ref {1.2e}). Since
$x^{(n,n) +1} = 0\in r_n(A)$, we have that $x$ is von Neumann
regular  by [Lemma 1]\cite {Ch81}.

(iii) Let $d_{ij}$ and $d_{ij}'$ in $A_{ij}$ denote the left and
right gm non-zero divisors of $A$ for any $i, j \in I$,
respectively. For any $x_{ij} \in D_{ij}$, there exists $u_{ts}
\in A_{ts}$ such that $d_{st}d_{ti} x_{ij}d_{js}'d_{st}'=
d_{st}d_{ti} x_{ij}d_{js}'d_{st}'u_{ts}d_{st}d_{ti}
x_{ij}d_{js}'d_{st}'$ and $ x_{ij}=
 x_{ij}d_{js}'d_{st}'u_{ts}d_{st}d_{ti} x_{ij}$
 since
$d_{st}d_{ti} x_{ij}d_{js}'d_{st}' \in B_{st}.$ This implies
$D_{ij} \subseteq r_n (A_{ij})$.
 Considering part (ii), we complete the proof.

(iv) If there exist $s, t \in I$ such that $r_n(A_{st}) \not=0$ ,
let $B_{st} = r _{n} (A_{st})$  and $D_{ij} = A_{is}B_{st}A_{tj} $
for any $i, j \in I.$ By part (iii), we have that $D=0$ and
$B_{st} =0.$ This is a contradiction. $\Box.$

\begin {Theorem} \label {1.3}   If $A$ has  left and
right gm non-zero divisors,  then $r_n(A)=g.m. r_n(A) = \sum \{r_n
(A_{ij}) \mid  i,  j \in I\}$ .
\end {Theorem}

{\bf Proof.} Let $B= r_n(A)$. For any $i, j \in I $ and $x_{ij}
\in B_{ij}$, there exists $y \in B$ such that $y _{ij} = x_{ij}.$
Let $d_{ii}$ and $d_{jj}'$ be left and right non-zero divisors in
$A_{ii}$ and $A_{jj}$, respectively. Since $(d_{ii} E(i,i))y
(d'_{jj}E(j,j)) = (d_{ii} x_{ij}d'_{jj})E(i, j) \in B,$ we have
that there exists $z \in B$ such that  $(d_{ii} x_{ij}d'_{jj})E(i,
j) = (d_{ii} x_{ij}d'_{jj})E(i, j) z (d_{ii} x_{ij}d'_{jj}) E(i,
j)$. By simple computation, we have   $d_{ii} x_{ij}d'_{jj} =
(d_{ii} x_{ij}d'_{jj}) z_{ji} (d_{ii} x_{ij}d'_{jj})$ and $ x_{ij}
= x_{ij}d'_{jj} z_{ji} d_{ii} x_{ij}$ . Thus $x_{ij}$ is von
Neumann regular . This implies $B_{ij} \subseteq r_n (A_{ij})$ and
$r_n (A) \subseteq \sum \{r_n (A_{ij}) \mid  i,  j \in I\}$.

Let $N = g.m.r_n(A).$ Since $g.m.r_n (A//N)=0$, we have that
$A_{ij}/N_{ij}$ is an  $r_n $-semisimple $A_{ji}/N_{ji}$-ring
for any $i, j \in I $ by Lemma \ref {1.2} (iv). It is clear that
$A_{ij}/N_{ij}$ is a $r_n $-semisimple $A_{ji}$-ring . This
implies $r_n (A_{ij}) \subseteq N_{ij}$ for any $i, j \in I.$
Consequently, $\sum \{r_n (A_{ij}) \mid  i,  j \in I\} \subseteq
g.m.r_n (A).$ $\Box$

If for any $s \in I$,  there exists $u_{ss} \in A_{ss}$ such
that $xe_{ss} = x$  for any $ i\in I$ and $x \in A_{is}$, then we
say that $A$ has a  right gm unit and $u_{ss}$ is  a right gm unit
in
 $A_{ss}.$ Similarly,
we can define a left gm unit of $A$  and gm unit of $A$. In
fact, if $A$ has  left and right gm units, then every ideal of $A$
is a gm ideal, so $r_n(A)= g.m.r_n(A) \subseteq \sum \{r
_n(A_{ij}) \mid  i,  j \in I\}$ by the proof of Theorem \ref
{1.3}.

It is clear that if $R$ is a ring and $M$ is a $\Gamma$- ring with
$R = M = \Gamma $, then $r_n (R) = r_n (M).$ We also have that $r
(R) = r (M)$   for $ r=r_b, r_k, r_l, r_j $ (see \cite [Theorem
5.2] {Bo87},
 \cite [Theorem 10.1] {CL71}  and \cite [Theorem 3.3] {Zh93}, \cite [Theorem 5.1] {CZ95}  ).

\begin {Theorem} \label {1.4} Let $r = r _b, r_l, r_j, r_n.$  Then

(i) $r(A)= g.m.r(A) = \sum \{r (A_{ij}) \mid  i,  j \in I\}$.

(ii) $r(A)=  \sum \{r (A_{ii}) \mid  i\in I\}$ when $A_{ij} = 0$
for any $i \not= j,$ i.e.  $r$ radical of the direct sum of rings
is equal to the direct sum of  $r$ radicals of these rings.

(iii) $r(A)$  is graded by $G$ when the index set $I$ of $A$ is an
abelian group $G$. Moreover  the grading is canonical.

Here $A$ has left and right  gm non-zero divisors when $r = r_n$
in  (i), (ii) and (iii).
\end {Theorem}

{\bf Proof.} (i)
\begin {eqnarray*}
r_b(A)&=& g.m.r_b(A) = \sum \{r_b (A_{ij}) \mid  i,  j \in I\}
\hbox { ( by \cite  [Theorem 3.7]{Zh93}) } \\
 r_l(A)&=& g.m.r_l(A) = \sum \{r_l (A_{ij}) \mid  i,  j \in I\} \hbox { \ \
 ( by \cite [Theorem 1.3 and Theorem 2.5]{CZ95} )}\\
r_j(A)&=& g.m.r_j(A) = \sum \{r_j (A_{ij}) \mid  i,  j \in I\}
\hbox { \ \ ( by  \cite [Theorem 3.10 and Theorem 1.3]{CZ95} ) }.
\end {eqnarray*}

(ii) Since the radicals of ring $A_{ii}$ and $A_{ii}$-ring
$A_{ii}$ are the same, we have (ii).

(iii) It follows from (i) and  \cite [Lemma 2.1] {ZZ03}. $\Box$

Let $M_I^f (R)$ denote the generalized matrix ring $A=\sum
\{A_{ij} \mid A_{ij} = R, i,  j\in  I\}$ with infinite index set $I$,
which is called an infinite matrix ring over ring $R$. In this
case,
 $M_I^f(k)$  is  called an infinite matrix algebra over field $k$. Let $M_{m \times n }(R)$ denote the ring
of all $(m \times n )$ matrices  over ring $R$.

\begin {Example}\label {1.5}
 (i)  Let $V = \sum _{g\in G}\oplus V_g$ be a
 vector space over field $k$ graded by abelian group $G$  with $dim \ V_g = n_g < \infty$.
Let $I$ denote  a basis of  $V.$
Then $\sum \{A_{ij} \mid   A_{ij} = Hom (V_j, V_i),   i,  j\in G\}
=\sum \{A_{ij} \mid  A_{ij} =  M_{n_i \times  n_j} (k), i, j\in
G\}  $ as generalized matrix  algebras. However,  $ \sum \{A_{ij}
\mid A_{ij} = Hom (V_j, V_i),   i,  j\in  G\} =\sum \{A_{ij} \mid
A_{ij} =  M_{n_i \times  n_j} (k), i,  j\in  G\}= \{ f \in End _k
V \mid
 ker f \hbox { has finite codimension } \}  = M_I^f (k)=
\sum \{A_{ij} \mid  A_{ij} = k, i,  j\in  I \} $ as algebras. Then
$r(M_I^f (k))= \sum \{ r(A_{ij}) \mid  A_{ij} =  k, i,  j\in I\}
 = r(\{ f \in End _k V \mid
 ker f \hbox { has finite codimension } \})=0 $ for $r = r_b, r _l , r_j$.
It is clear that generalized matrix algebra $A= \sum \{A_{ij} \mid
A_{ij} = Hom (V_j, V_i),   i,  j\in  G\}$ has left and right gm
non-zero divisors if and only if $n_i = n_j $ for any $i, j \in
G.$ Consequently, $r_n( M_I^f (k))= \sum \{r_n(A_{ij}) \mid A_{ij}
= k, i,  j\in  I \} = M_I^f(k)$. That is, $\{ f \in End _k V \mid
 ker f \hbox { has finite codimension } \}$ is a von Neumann regular algebra.

(ii) By Theorem \ref {1.4}, $r(M_I^f(R)) = M_I^f (r(R))$ for $r=
r_b, r_l, r_j.$ If $R$ is  a ring with left and right non-zero
divisors, then  $r_n(M_I^f(R)) = M_I^f (r_n(R)).$ Obviously,  if
$R$ has left and right units,
 then $R$ is a ring with left
and right non-zero divisors ($R \not= 0$), so $r_n(M_I^f(R)) =
M_I^f (r_n(R)).$
\end {Example}

\section {Application in path  algebras  }

\begin {Lemma} \label {3.1}  Let $r$
denote $r_b ,  r_k ,  r_1, r_j$ and $s, t \in I$ . If $A_{st}
\not= 0$, then

(i)  $r(A_{st}) = 0$ if and only if $A_{ts}\not= 0$.

(ii)  $r(A_{st}) = A_{st}$ if and only if $A_{ts}=0$.
\end {Lemma}
{Proof.} If $r (A_{st})= 0$, then $r _b(A_{st})=0 $ and $A_{ts}
\not=0.$ Conversely, if $A_{ts}\not= 0$ and $r_j (A_{st}) \not=0,$ then $r_j (A_{st}) = k$ or
 there exists $y \in r_j(A_{st})$ with $l (y) > 0$. Since $y$ is a right
quasi-regular element of $A_{ts}$-ring $A_{st}$,   for
$x\in A_{ts}$,   there
exists $u \in A_{ts}A_{st}$ such that
\begin {eqnarray} \label {3.1e} y(xy)u + y(xy) &=&  -yu.
\end {eqnarray}
If $l(u)>0$, then  the right hand side of (\ref {3.1e}) is shorter
than the left hand side of (\ref {3.1e}), we get a contradiction.
If $l(u)=0$, then  the left hand side of (\ref {3.1e}) is either
equal to zero, or longer than the right hand side of (\ref
{3.1e}).  We get a contradiction. This implies that $r_j
(A_{st})=0$ and $r(A_{st})= 0 .$ $\Box$

\begin {Lemma} \label {3.2}
$ r_n (A_{st}) = 0$ for any $s\not= t.$
\end {Lemma}
{Proof.} For any $0 \not=x \in A_{st}$,  if $x$ is a von Neumann
regular element, then there exists $y\in A_{ts}$ such that $x=
xyx$. Considering the length  of both sides we get a
contradiction. Consequently, $r_n (A_{st})=0.$ $\Box$

 \begin {Theorem} \label {3.3} (i)
 $r (A)=g.m. r(A)= \sum \{r(A_{ij}) \mid i, j\in I \} = k R(D)$,
\noindent where $r$ denotes $r_b, r_1,r_k$ and  $r_j$.

(ii)  $r_n (A)=g.m. r_n(A) = \oplus  \{N_{ii} \mid i\in I \},$\\
where $ N_{ii}= \left \{
\begin {array} {ll} ke_{ii}= r_n(A_{ii})
  & \ \ \hbox { when }  A_{si}= A_{is} =0 \hbox { for any  }  s \in I  \hbox { with  } i \not= s. \\
  0 & \ \ \hbox {otherwise   }.
  \end {array} \right. $\\
\end {Theorem} {\bf Proof.}  (i) By Lemma \ref {3.1},  $\sum \{ r(A_{ij})
\mid i, j\in I \} \subseteq  k R(D)$. Let $x \in k R(D)$ be a
regular path from $i$ to $j$. Then $A_{ij} \not=0$ and $A_{ji}=0,
$ which implies $r(A_{ij})= A_{ij}$
 and $x \in r(A_{ij}.$ This has proved $\sum \{r(A_{ij}) \mid i, j\in I \} = k R(D)$. It follows that
$r_b (A)= r_j(A)= k R(D)$ from \cite [Proposition 5] {Lu88} or
Theorem \ref {1.4}. Thus $r (A)=  k R(D)$. By \cite [Theorem 1.3]
{CZ95} or Theorem \ref {1.4}, $r(A) = g.m.r(A).$ We complete the
proof.

(ii) Since $A$ has a gm unit $e_{ss}\in A_{ss}$ for any $s\in I$,
we have that $r_n (A) = g.m. r_n(A)$.
 By Lemma \ref {3.2} and the proof of Theorem \ref {1.3}, we have $r_n (A)=g.m. r_n(A)
\subseteq  \sum \{r_n(A_{ij}) \mid i, j\in I \}  = \oplus
\{r_n(A_{ii}) \mid i\in I \}$. Let $N = r_n (A)$. If  $N_{ss}
\not=0$, then  $ N_{ss} = r_n (A_{ss}) = ke_{ss}$ by the proof of
Lemma \ref {3.2}. For any $t \in I$ with $t \not= s$, since
$A_{ts} N_{ss} = A_{ts} \subseteq N_{ts} = 0$, we have $A_{ts}=
0.$ Similarly, $A_{st} =0.$
 $\Box$

Next we give the  relations  between the radicals of path algebras and
connectivity of directed graphs.

\begin {Theorem} \label {3.4} Directed graph $D$ is strong connected if and only if
$A=A(D)$ is a prime algebra.
\end {Theorem}
{\bf Proof.} If $D$ is strong connected  and $A$ is not prime,
then there exist $u, v , s, t \in I$ and $0 \not= x \in A_{uv},
0\not= y \in A_{st}$ such that $x Ay =0$, i.e. $x A_{vs}y =0$,
which contradicts with the strong connectivity of $D$. Consequently $A$ is
prime. Conversely, if $A$ is prime, then $e_{ii}A_{ij} e_{jj}
\not=0$ for any $i, j \in I,$ which implies that $D$ is strong
connected. $\Box$

\begin {Theorem} \label {3.5} Every weak component of  $D$  has  at least two vertexes if and only if
 $r_n (A)=0$.
\end {Theorem}

{\bf Proof.} It follows from Theorem \ref {3.3} (ii). $\Box$

\begin {Lemma} \label {3.6}  Every  directed graph  $D$  is the union of all of its  unilateral
components.
 \end  {Lemma}
{\bf Proof.} For any path $x \in A_{st}$, set $${\cal K} = \{ E
\mid E  \hbox { is a unilateral connected subgraph of $D$ with }
x \in E\}.$$ By Zorn's Lemma, we have that there exists a maximal
$Q$ in ${\cal K}$. $\Box$
 \begin {Theorem} \label {3.7}  Let $r$ denote $r_b,  r_k$, $r_l$ and $r_j$,  respectively.
The following conditions are equivalent.

(i) Every weak component of $D$ is a strong component.

(ii) Every unilateral component of $D$ is a strong component.

(iii) Weak component,  unilateral component  and  strong component
of $D$ are the same.

(iv) $D$ is the union of strong components  of $D$.

(v) $D$ has no any regular path.

(vi) $A_{ij} = 0$ if and only if $A_{ji} = 0 $ for any $i, j \in
I.$

(vii) $A$ is a direct sum of prime  algebras.

(viii) $A$ is semiprime.

(ix)
 $A_{ij}$ is a semiprime $A_{ji}$-ring for $i, j\in I$.

(x) $r(A_{ij}) = 0$ for any $i, j \in I$.

(xi) $r(A)= 0.$

\end {Theorem}
{\bf Proof.} By Theorem \ref {3.3}, (v), (vi), (viii), (ix), (x)
and (xi) are equivalent.

$(i) \Rightarrow (vi)$. If  $i$ and $j$ belong  to the same weak
component, then $A_{ij} \not=0$ and $A_{ji} \not= 0$. If  $i$ and
$j$ do not belong to  the same weak component, then  obviously
$A_{ij} =0$ and $A_{ji} = 0$.

$(ii) \Rightarrow (vi)$. If $A_{ij}\not= 0,$ then  there exists a
path  $x \in A_{ij}.$ By Lemma \ref {3.6}, $x$ belongs to a certain
unilateral component of $D$. Consequently,  $x$ belongs to a certain
strong component of $D$. This implies $A_{ji} \not=0. $

$(vi) \Rightarrow (ii)$. If $i$ and $j$ belong to the same unilateral
component of $D$, then $A_{ij} \not= 0$
 or $A_{ji} \not=0$. Consequently $A_{ij} \not=0 $ and $A_{ji} \not=0$, which  implies that  $i$ and $j$  belong to the same strong component of $D$.
Therefore (ii) holds.

$(iv) \Rightarrow (vi)$. If $A_{ij} \not= 0$, then there exists a
path $x \in A_{ij}$ and  $x$ belongs to a certain  strong component
of $D$. This implies $A_{ji} \not=0.$

$(vi) \Rightarrow (iv)$.  For any arrow  $x \in A_{ij}$, we only
need show that
 $x$ belongs to a
certain  strong component of $D$. By Lemma \ref {3.6}, there
exists a certain unilateral component $C$   of $D$ such that $x
\in C$. Since (ii) and (vi) are equivalent, we have that $C$ is a
strong component of $D$.

$(iv) \Rightarrow (i)$. If $i$ and $j$ belong to the same weak
component, then there exists a semi-path $x = x_{i i _1}
x_{i_1i_2} \cdots x_{i_nj}$. If $i$ and $j$ belong to different
strong components, then we can assume that $i_s$ is the first
vertex, which does not belong to the strong component containing
$i$. Consequently,
   $A_{i_{s-1}i} \not=0$  and    $A_{i i_{s-1}} \not=0$, and   either  $A_{ii_s } \not=0$ or
$A_{i_si} \not=0$. Since  (iv)  and  (vi) are equivalent, we have that
$A_{ii_s } \not=0$ and    $A_{i_si} \not=0$. We get a
contradiction. This shows that  $i$ and $j$ belong to the same
strong components.

 $(iii) \Rightarrow (i)$ is obvious.

$(i) \Rightarrow (iii)$. Since (i) and (ii) are equivalent, we
have (iii).

    $(vii) \Rightarrow (viii)$. It follows from Theorem \ref {1.4} (ii).

$(iv) \Rightarrow (vii)$. Let $\{ D^{( \alpha )} \mid \alpha \in
\Omega \}$ be all of the strong component of $D$ and $D =\cup \{
D^{(\alpha )} \mid \alpha \in \Omega   \}$. Thus $A(D) = \oplus
\{ A( D^ {(\alpha )} ) \mid \alpha \in \Omega   \}$. However, for
any $\alpha \in \Omega $, $A(D^ {(\alpha )} )$ is a prime algebra
by Theorem \ref {3.4}. We complete the proof. $\Box$

We easily obtain the following  by the preceding conclusion for $
r= r_b, r_l, r_k, r_j$.
 $D$ has no cycle  if and only if
$r(A_{ij})=A_{ij}$ for any  $i\not=  j\in  I$;
  $s$ and $t$ ($s \not= t$) are not contained in the same
cycle if and only if $r(A_{st})=A_{st}$; $s$ and $t$ ( $s \not=
t$ ) are contained in the same  cycle if and only if
$r(A_{st})=0$.

We give an example to show whether the condition in Theorem \ref
{1.3} is
 a necessary one.
\begin {Example}\label {3.8} (i)
Let $D$ be a directed graph with vertex set  $I = \{1, 2\}$ and
only one arrow $x_{12} \in A_{12}.$ Obviously, $A_{12} =k x_{12},
A_{11} =k e_{11}, A_{22} = k e_{22}, A_{21} =0,  $
 $r_n(A_{ii}) = ke _{ii}$ and    $r_n(A_{ij}) = 0$  for any $i, j \in I$ with $i \not= j.$
By Theorem \ref {3.3} (ii), $r_n (A) = 0 \not= \sum \{ r_n
(A_{ij}) \mid i, j \in I \}.$ It is clear that  $A$ has no left gm
non-zero divisor since $A_{21} =0.$ Consequently, it is possible
that Theorem \ref {1.3} does not hold if its condition is dropped.

(ii) Let $I = \{1, 2\}$ and $A_{ij} = M_{i \times j} (k)$ for any
$i, j \in I.$ It is clear that $A$ has no left gm non-zero divisor
in $A_{12}$ since, for any non-zero $x\in A_{12}$, there exists a
non-zero $y\in A_{21}$ such that $xy = 0 $. However, $r_n (A) =
\sum \{ r_n (A_{ij}) \mid i, j \in I \} = M_{3 \times 3} (k).$
Consequently,  the condition in Theorem \ref {1.3} is not a
necessary condition.
\end {Example}

\section {Application in networks}

 We give an example to establish that  (generalized) path algebras are
applied in networks.
\begin {Example} \label {6.1}
Let
\begin {eqnarray*}I_N &=& \{x \mid x \hbox { is the name of a country }\},\\
I_O &=& \{x \mid x \hbox { is the name of an organization} \},\\
I_S &=& \{x \mid x \hbox { is the name of a server } \},\hbox {\ \ \ and}\\
I_C &=& \{x \mid x \hbox { is the name of a computer } \}.\\
\end {eqnarray*} Let $I$ denote one of $I_N, I_O, I_S$ and  $I_C$. We construct a directed graph $D$
as follows. Let $I$ be the vertex set of $D$ and
each  piece of  information from $i$ to $j$  denote an arrow from $i$ to $j$ for $i, j\in I$.
We get a path algebra.
\end {Example}

By the way we give an example of  applications of braided diagrams.
\begin {Example} \label {6.2}
Let $V_1,V_2,V_3 $ and  $ V_4 $ be  vector spaces with basis
\begin {eqnarray*}& & \{x \mid x \hbox { is the number of a train }\},\\
 & & \{x \mid x \hbox { is the number of an automobile} \},\\
 & & \{x \mid x \hbox { is the number of a boat } \} \hbox {\ \ \ and} \\
 & & \{x \mid x \hbox { is the number of an airplane } \},\hbox {respectively}.\\
\end {eqnarray*} Let $G= {\bf Z}_4$ and $r$ a bicharacher on $G$. In this way, we can turn the
map about pathes of trains, automobiles, boats and airplanes into a braided diagram. Note the following
details:

(i) If the path is parallel  to the latitude, then we view that the west is higher than the east.

(ii) A braiding denotes a cross of two pathes.

(iii) If the path is concave, then we view that the path is evaluation.

(iv) If the path is convex, then we view that the path is
coevaluation.
\end {Example}

\vskip 0.3cm {\bf Acknowledgement }: The work is supported by
Australian Research Council. S.Zhang. thanks the Department of
Mathematics, University of Queensland for hospitality.

\end {document}